%% file: sch_arxiv_v2.tex
\documentclass[11pt,reqno]{amsart}

\usepackage{mathtools}
\numberwithin{equation}{section}
\usepackage{amscd,amssymb,amsthm,amsbsy}
\usepackage{thmtools}
\usepackage[mathscr]{euscript}
\usepackage{dsfont}

\usepackage[width=5.65in,height=8.0in,centering]{geometry}
\usepackage{xcolor}
\usepackage{booktabs}
\usepackage{microtype}
\usepackage[shortlabels]{enumitem}
\usepackage{comment}
\usepackage{subcaption} 

\usepackage{tikz}
\usetikzlibrary{calc,intersections,arrows,patterns,cd,shapes.misc,arrows.meta}

\usepackage[
	colorlinks=true,
	linkcolor=blue,
	citecolor=purple,
	urlcolor=blue,
	linktoc=page,
	plainpages=false,
	pdfpagelabels,
	backref=page
]{hyperref}
\usepackage[nameinlink,noabbrev]{cleveref}


\theoremstyle{plain}
\newtheorem{theorem}{Theorem}[section]
\newtheorem{lemma}[theorem]{Lemma}
\newtheorem{corollary}[theorem]{Corollary}

\newtheorem{proposition}[theorem]{Proposition}
\newtheorem{propdef}[theorem]{Proposition/Definition}

\newtheorem{question}[theorem]{Question}

\theoremstyle{definition}
\newtheorem{definition}[theorem]{Definition}
\newtheorem{notation}[theorem]{Notation}
\newtheorem{exm}[theorem]{Example}

\theoremstyle{remark}
\newtheorem{rem}[theorem]{Remark}

\newenvironment{example}%
{\pushQED{\qed}\begin{exm}}%
{\popQED\end{exm}}  

\newenvironment{remark}%
{\pushQED{\qed}\begin{rem}}%
{\popQED\end{rem}}  

\newcommand{\abs}[1]{\left|#1\right|}     
\newcommand{\set}[1]{\left\{#1\right\}}   
\newcommand{\angl}[1]{\left\langle #1 \right\rangle}    

\newcommand{\newword}[1]{\emph{#1}}

\DeclareMathOperator{\OS}{{OS}}       
\DeclareMathOperator{\IN}{{IN}}       
\DeclareMathOperator{\IA}{{IA}}       
\DeclareMathOperator{\EA}{{EA}}       
\DeclareMathOperator{\In}{{In}}       
\DeclareMathOperator{\BC}{{BC}}       
\DeclareMathOperator{\bnbc}{{\beta{\mathbf{nbc}}}}  
\DeclareMathOperator{\GL}{{GL}}       
\DeclareMathOperator{\bc}{{\overline{BC}}}       
\DeclareMathOperator{\os}{{\overline{OS}}}    

\DeclareMathOperator{\sh}{{sh}}       
\DeclareMathOperator{\Sh}{{Sh}}       
\DeclareMathOperator{\link}{{link}}   

\newcommand{\irr}{\mathrm{irr}}     
\newcommand{\rS}{\overline{S}}      
\newcommand{\rE}{\overline{E}}      

\def\bfcdot{{\hbox{$\scriptscriptstyle\bullet$}}} 

\newcommand{\CC}{\mathds{C}}     
\renewcommand{\k}{\Bbbk}
\newcommand{\RR}{\mathbb{R}}     
\newcommand{\ZZ}{{\mathbb{Z}}}   
\newcommand{\EE}{{\bigwedge}}    
\newcommand{\OLambda}{{\overline{\Lambda}}} 

\newcommand{\A}{{\mathscr{A}}}  
\newcommand{\B}{\mathscr{B}}    
\newcommand{\E}{\mathscr{E}}    
\newcommand{\Z}{\mathscr{Z}}    
\newcommand{\sS}{\mathscr{S}}    
\renewcommand{\L}{\mathscr{L}}  
\newcommand{\M}{\mathsf{M}}     
\newcommand{\U}{\mathsf{U}}     
\newcommand{\R}{\mathsf{R}}	
\newcommand{\nbc}{\mathbf{nbc}}

\newcommand{\bB}{\mathsf{B}} 

\renewcommand{\P}{\mathbb{P}}
\newcommand{\C}{\mathbb{C}}

\DeclareMathOperator{\cl}{{cl}}  
\DeclareMathOperator{\rank}{{rk}}  
\DeclareMathOperator{\nullity}{{null}}  
\DeclareMathOperator{\ann}{{ann}}  
\DeclareMathOperator{\Hom}{Hom}

\newcommand{\uppercaseFactory}[2]{
    \foreach \x in {A,...,Z}{
        \expandafter\xdef\csname #1\x\endcsname{#2{\x}}}}
\newcommand{\lowercaseFactory}[2]{
    \foreach \x in {a,...,z}{
        \expandafter\xdef\csname #1\x\endcsname{#2{\x}}}}

\uppercaseFactory{sh}{\noexpand\mathcal}
\uppercaseFactory{cal}{\noexpand\mathcal}
\uppercaseFactory{scr}{\noexpand\mathscr}
\uppercaseFactory{mod}{}
\uppercaseFactory{fr}{\noexpand\mathfrak}
\lowercaseFactory{fr}{\noexpand\mathfrak}
\uppercaseFactory{rm}{\noexpand\mathrm}
\lowercaseFactory{rm}{\noexpand\mathrm}


\begin{document}
\begin{abstract}
  We compute the resonance varieties of the Orlik--Solomon algebras of
  Schubert matroids.  Previously known results have been largely confined to
  matroids of rank $3$, and the Schubert matroids give a family of
  nontrivial examples of arbitrary rank.  We give an expression for the
  cohomology of the Orlik--Solomon algebra that we show
  holds for Schubert matroids
  as well as for simple, complex-realizable rank-$3$ matroids that do not
  support multinets.
  Our proof introduces a family of matroids for which the 
  Orlik--Solomon algebra is isomorphic to an exterior Stanley--Reisner algebra.
\end{abstract}
  
\title[Schubert arrangements]{Resonance varieties of Schubert arrangements}
\author{Graham Denham, Prajwal Udanshive}
\maketitle
\tableofcontents

\section{Introduction}\label{sec:intro}
Let $\M$ be a matroid on an ordered ground set $E$, and let $\k$ be a field of
characteristic different from $2$.  The 
Orlik--Solomon algebra of $\M$, denoted $\OS(\M)$, is a finite-dimensional
graded $\k$-algebra whose Hilbert series equals (up to signs and reindexing)
the characteristic polynomial of $\M$.  Originally, $\OS(\M)$ was introduced
for complex hyperplane arrangements as a presentation for the cohomology of
the arrangement complement~\cite[Thm.\ 5.2]{orlik-solomon}.
From a more recent point of view, it has been studied purely as a matroid
invariant.  For example, 
it is known to be valuative, in a suitably generalized sense~\cite{EMPV25},
and Coron~\cite{Co25} showed it carries the structure of a Feynman operad.

Here, we are interested in resonance for matroids.  In the motivating
situation where $\M=\M(\A)$
is the matroid of a complex hyperplane arrangement $\A$,
one is interested in computing $H^\bfcdot(U(\A),\CC_\rho)$
in terms of $\A$ and $\CC_\rho$, a rank-$1$ local system on the arrangement complement
$U(\A)$ (see~\cite{ESV92,STV95}).  We refer to Suciu's surveys~\cite{Su16,Su25}
for a full discussion.  Cohomology of local
systems $\CC_\rho$ in the neighbourhood of the trivial representation
$\CC$ can be approximated by a calculation using 
only the cohomology ring $\OS(\M)=H^\bfcdot(U(\A),\CC)$.

More precisely, the image of any $w\in\k^E$ is a degree-$1$ element of
$\OS(\M)$.  Graded-commutativity gives $w^2=0$, so the
map $x\mapsto w\cdot x$ endows $\OS(\M)$ with the structure of a cochain
complex.  The cohomology of this complex 
is a matroid invariant;
by comparison, it is unknown in general if $H^\bfcdot(U(\A),\CC_\rho)$ is
determined by $\M(\A)$.

For generic $\rho$ and $w$, the cohomology of local systems $\CC_\rho$
and of the Orlik--Solomon algebra vanishes, 
so it makes sense to study the non-vanishing loci: these are the
\emph{characteristic varieties} and \emph{resonance varieties}, respectively,
of the arrangement and matroid.
A characterization of $\M$ and $w$ for which $H^1(\OS(\M),w)\neq 0$ is
given by Falk and Yuzvinsky~\cite{FY07}, in terms of the combinatorics of
multinets, at least for complex-realizable matroids.  They reduce to the
rank-$3$ case; however, their
characterization is somewhat involved.  The resonance varieties of
matroids are not well-understood beyond cohomological degree $1$, and
we refer to the surveys~\cite{Su25,De16} for more discussion.

Resonance for uniform
matroids $\U_{r,n}$ was computed by Yuzvinsky~\cite{Yu95} for all
$1\leq r\leq n$: cohomology
is concentrated in degree $r-1$.  He showed
\begin{equation}\label{eq:uniform}
H^\bfcdot(\os(\U_{r,n}),w)\cong \k^{\beta(\U_{r,n})}[1-r]
\end{equation}
as a graded vector space, for all nonzero $w$.  (Here, $\os$ denotes the reduced
Orlik--Solomon algebra.)

The objective of this paper is to give a computation of resonance in
all cohomological degrees for the Schubert matroids.  These matroids
are structurally
straightforward; however, they play a key role in the study of
matroid valuations~\cite{derksenfink10}.  Schubert matroids (properly) include 
uniform matroids, so our results generalize \eqref{eq:uniform}.
The formula~\eqref{eq:res formula} we obtain is also satisfied by some
matroids that are not Schubert matroids, which we discuss in \S\ref{sec:last}.

\begin{definition}\label{def:res formula}
  We say that a matroid $\M$ satisfies the \newword{convolution resonance formula}
  over $\k$ if, for all nonzero $w\in\Delta_E$, 
  we have
    \begin{equation}\label{eq:res formula}
  H^\bfcdot(\os(\M),w) \cong \bigoplus_{\substack{
  F\in \L^\irr(\M)\colon\\
  w\in \Delta_F}}
  \k^{\beta(\M\vert_F)}\otimes_{\k} \OS(\M/F)[1-\rank_\M(F)].
    \end{equation}
Here, the sum runs over connected flats $F$ of $\M$, and $\Delta_F$ is
a linear subspace of $\os^1(\M)$ of dimension $\abs{F}-1$ which we define
below.
\end{definition}

\begin{theorem}\label{thm:schres}
  Loopless Schubert matroids satisfy the convolution resonance
  formula~\eqref{eq:res formula}.
\end{theorem}

Our proof makes use of an observation that Schubert matroids have particularly
simple Orlik--Solomon algebras.  We introduce the class of \emph{arboreal}
matroids, which are defined by the property that their positive-rank,
connected flats are \emph{tree-convex}: that is, there exists a
tree with vertices
$E$ for which each such flat is (the vertex set of) an induced subtree.

Tree-convex sets were independently characterized by Duchet~\cite{Du78},
Flament~\cite{Fl78} and Slater~\cite{Sl78} and revisited recently by
Chudnovsky--Nguyen--Scott--Seymour~\cite{chud26}.
A collection of subsets $\sS$ has the
\emph{Helly property} if, whenever $F_1,\ldots,F_k\in \sS$ and
$F_i\cap F_j\neq\varnothing$ for all $1\leq i<j\leq k$, then
$F_1\cap \cdots\cap F_k\neq \varnothing$: i.e., the nerve of $\sS$ is a
clique complex.  Let
$\L^\irr_{\geq1}(\M)$ denote the connected flats of $\M$ of positive rank.
Their characterization is the following.
\begin{theorem}\label{thm:arboreal-chud26}
  A matroid $\M$ is arboreal if and only if the intersection graph of
  $\L^\irr_{\geq1}(\M)$ is chordal, and $\L^\irr_{\geq1}(\M)$ has the Helly property.
\end{theorem}

Examples include Schubert matroids, laminar matroids, and free coextensions.
For any matroid, it is well-known that the Orlik--Solomon algebra $\OS(\M)$
has a Gr\"obner degeneration to the exterior face ring of the
broken-circuit complex of $\M$ (see \S\ref{subsec:BC}).  Arboreal matroids
have the special property that this degeneration is actually an
isomorphism:
\begin{theorem}\label{thm:degeneration}
  If $\M$ is an arboreal matroid,
  \[
  \OS(\M)\cong\k\angl{\BC(\M)},
  \]
  with respect to any compatible ordering\footnote{to be made precise in Section~\ref{ss:arboreal}.} 
  of the ground set $E$.
\end{theorem}
As a result, we see that all matroid exterior face rings are
Orlik--Solomon algebras.  This ``categorifies'' a result of 
Brylawski~\cite[Thm.\ 4.2]{Bry77} who showed that, for any matroid $\M$, the
independence complex $\IN(\M)$ is isomorphic to the (reduced) broken circuit
complex of the free coextension $0\times_E M$.
\begin{corollary}
  For any matroid $\M$, we have $\os(0\times_E \M)=\k\angl{\IN(\M)}.$
\end{corollary}

\section{Orlik--Solomon algebras}\label{sec:setup}
Write $[a,b]:= \set{a, a+1, \dots, b}$ for integers $a\leq b$.
Throughout, we fix a field $\k$ of characteristic not $2$.
Our matroids
have (finite) ordered, nonempty ground sets
$E$ which we will often identify with the set
$[n]:=\set{0,1,2,\ldots,n}$. 
Let $\set{e_i\colon i\in E}$ denote a basis for the vector space $\k^E$, and
$\set{x_i\colon i\in E}$ the dual basis.  

\subsection{Definitions}
\begin{definition}\label{def:Delta}
Let $\partial\colon \k^E\to\k$
denote the linear map given by $\partial(e_i)=1$ for each $i\in E$, and 
let $\Delta_E:=\ker\partial$.
Clearly, $\Delta_E$ is given by the equation $\sum_{i\in E}x_i=0$.

More generally, for any subset $S\subseteq E$, we let
\begin{equation}\label{eq:coordinate-subspace}
\Delta_S:=\set{x\in \k^E\colon \sum_{i\in S}x_i=0,\text{~and~
$x_i=0$ for $i\in E\setminus S$}}.
\end{equation}
\end{definition}

The linear map $\partial$ extends to a derivation of degree $-1$ on the exterior
algebra $\Lambda_E:=\EE(\k^E)$.  We let $\OLambda_E$ denote the kernel of the
derivation, which is easily seen to be the exterior subalgebra,
$\OLambda_E=\EE(\Delta_E)$.

For a subset $S\subseteq E$, we abbreviate
\[
e_S := e_{s_1}e_{s_2}\cdots e_{s_k}
\]
where $S=\set{s_1<\cdots<s_k}$.

\begin{lemma}\label{lem:partial_product}
  For any nonempty $S\subseteq E$, let $T$ be a directed tree with vertices
  $S$.  The set
\[
\bB(T):=\set{e_i-e_j\colon (i,j)\in\E(T)}
\]
  forms a basis of $\Delta_S$, and 
\begin{equation}\label{eq:partial_product}
\partial(e_S) = \pm \prod_{(i,j)\in\E(T)}(e_i-e_j),
\end{equation}
where $\E(T)$ denotes the set of edges of $T$.
\end{lemma}
\begin{proof}
  Since each $e_i-e_j\in\Delta_S$ and $\dim\Delta_S=\abs{S}-1$, the first
  assertion follows by noting
  $\bB(T)$ is linearly independent, by the acyclicity of $T$.
  Then the right-hand side of \eqref{eq:partial_product}
  spans $\bigwedge^{\abs{S}-1}\Delta_S\cong \k$.
  Noting that monomials on both sides appear with coefficients $\pm 1$, it
  is enough to check $\partial(e_S)\in\bigwedge^{\abs{S}-1}\Delta_S$, which follows from
  an elementary verification that $(e_i-e_j)\partial(e_S)=0$ for all
  $e_i-e_j\in \Delta_S$.
\end{proof}

\begin{definition}\label{def:red_OS}
For a matroid $\M$ and field $\k$, the \newword{reduced Orlik--Solomon algebra}
$\os(\M)$ is the graded $\k$-algebra
\[
\os(\M):=\OLambda_E/I_{\M},
\]
where
\[
I_{\M}:=\angl{\partial(e_C)\colon \text{$C$ is a circuit of }\M}
\]
is called the \newword{Orlik--Solomon ideal}. The \newword{Orlik--Solomon
algebra} $\OS(\M)$ is the extension of $\os(\M)$ to $\Lambda_E$.  For each
$i\in E$, the \newword{affine Orlik--Solomon algebra} $\OS(\M;i)$ is the quotient
$\OS(\M)/\angl{e_i}$.
\end{definition}
Here are some classical facts about Orlik--Solomon algebras, for which we
  refer to \cite{OTbook}.
\begin{proposition}\label{prop:loops}~
  \begin{enumerate}
\item
  Each affine Orlik--Solomon algebra
  admits an isomorphism $\phi_i\colon \OS(\M;i)\to\os(\M)$ given in degree $1$ by
  $\phi_i(e_j)=e_j-e_i$ for all $j\in E\setminus\set{i}$.
\item
  If $i\in E$ is a loop, then $\set{i}$ is a circuit, which implies
  $1\in I_{\M}$, and $\OS(\M)=\os(\M)=0$.  For this reason, we restrict our
  attention to loopless matroids for the rest of the paper.
\item
  If $i,j\in E$ are parallel, then $\set{i,j}$ is a circuit, which implies
  $e_i=e_j$ in $\OS(\M)$.  Up to isomorphism, then, it is enough to
  consider simple matroids.
\item
  The derivation $\partial$ is well-defined on $\OS(\M)$, so
  $(\OS(\M),\partial)$ is
  a chain complex of $\k$-vector spaces.  It is split exact, and
  $\ker\partial\cong\os(\M)$.
  \end{enumerate}
  \end{proposition}
  \begin{remark}
These combinatorial facts reflect the topological origins.
  If $\M$ is the matroid of a (central) complex hyperplane arrangement in
  $\C^r$, then $\OS(\M)$ is isomorphic to the cohomology algebra
  of the complement of the union of hyperplanes in $\C^r$,
  while $\os(\M)$ and $\OS(\M;i)$ are the cohomology algebra of the complement
  in $\P^{r-1}$, or equivalently  
  in $\C^{r-1}=\P^{r-1}\setminus{H_i}$, where $H_i$ denotes the $i$th
  hyperplane.
  The derivation $\partial$ is given by contraction along the Euler vector
  field: for details, we refer to \cite[\S6.1]{Dimcabook}.
\end{remark}
\subsection{Alternative presentations}
The Orlik--Solomon algebra was originally defined in terms of circuits.  It
will be convenient for us to use the following cryptomorphic definition, given
in terms of connected flats. 
\begin{lemma}\label{lem:flat}
  Let $\M$ be a matroid on $E$.  For a nonempty subset $S\subseteq E$,
  the ideal $\angl{\Delta_S}$ satisfies
\[
\angl{\Delta_{S}}^{\rank(S)}\subseteq I_{\M}.
\]
\end{lemma}
\begin{proof}
  Let $d=\rank(S)$.  Since the ideal $\angl{\Delta_S}$ is
  generated in degree $1$ by differences $e_s-e_{s'}$ for $s,s'\in S$,
  its power $\angl{\Delta_S}^d$ is generated by elements $\partial(e_T)$ for
  sets $T\subseteq S$ with $\abs{T}=d+1$, by Lemma~\ref{lem:partial_product}.
  Since $\rank(T)\leq \rank(S)=d$, the set $T$ is dependent in $\M$, so it
  contains a circuit $C$.  It follows that $\partial(e_T)$ is a multiple of
  $\partial(e_C)$, which implies $\partial(e_T)\in I_{\M}$.
\end{proof}

\begin{theorem}\label{thm:OS_pres}
Let $\M$ be a matroid.  Let $\sS$ be any collection of
nonempty subsets of $E$ that contains all connected flats. Then
\begin{equation}\label{eq:OS_pres}
\os(\M)=\OLambda_E/\big(\sum_{F\in\sS}\angl{\Delta_F}^{\rank_{\M}(F)}\big).
\end{equation}
\end{theorem}
\begin{proof}
  If $C$ is a circuit of $\M$, its closure $F=\cl(C)$ is a connected flat of
  rank $\abs{C}-1$ (e.g.,~\cite[Lem.\ 2.8]{BdMN03}), and we see from
  \eqref{eq:partial_product} that
  $\partial(e_C)\in\angl{\Delta_F}^{\abs{C}-1}$.  Since $\sS$ contains all
  connected flats, 
\[
I_{\M}\subseteq \sum_{F\in\sS}\angl{\Delta_F}^{\rank_\M(F)}.
\]
Conversely, applying Lemma~\ref{lem:flat} to each nonempty set $F\in \sS$,
we find
\[
\sum_{F\in\sS}\angl{\Delta_F}^{\rank_\M(F)} \subseteq I_{\M}.
\]
\end{proof}
For a graded $\k$-linear space $V$, we grade its dual by
$(\Hom_{\k}(V,\k))_p:= \Hom_\k(V_{-p},\k)$ for all $p\in\ZZ$.
\begin{corollary}\label{cor:OS_dual}
  For any collection of subsets $\sS$ containing the connected flats, we
  have an isomorphism of graded $\Lambda(\Delta_E)$-modules
  \[
  \Hom_{\k}(\os(\M),\k)[1-|E|]\cong \ann(I_{\M})=\bigcap_{F\in\sS}\angl{\Delta_F}^{\nullity_\M(F)}.
  \]
\end{corollary}

\begin{proof}
The first isomorphism holds for every homogeneous ideal $I\subseteq
\OLambda_E$ because $\OLambda_E$ is a Frobenius algebra. Let $W$ be a linear
subspace of $\Delta_E$, and $L:=\angl{W}$ be the associated linear
ideal in $\OLambda_E$. Recall the annihilator of $L^k$ is
$L^{\dim(W)-k+1}$. Thus, each
\begin{align*}
  \ann(\Delta_F^{\rank(F)}) &= \angl{\Delta_F}^{\abs{F}-1-\rank(F)+1}\\
  &=\angl{\Delta_F}^{\nullity_\M(F)}.
\end{align*}
The conclusion follows from Theorem~\ref{thm:OS_pres}.
\end{proof}

For the sake of completeness, we remark that the Orlik--Solomon algebra also
has a presentation via matroid bases
as an intersection of ideals generated in degree $1$.
If $B=\set{b_1,\ldots,b_r}$ is a basis and
$b_1<b_2<\cdots<b_r$, for $0\leq i\leq r$, let
\[
F(B)_i=b_1\vee b_2 \vee \cdots \vee b_i,
\]
so that $\varnothing=F(B)_0$ and $E=F(B)_r$.  Further, let
  \[
  \Delta_{F(B)}\coloneqq \sum_{i=1}^r\Delta_{F(B)_i\setminus F(B)_{i-1}}
  \subseteq \Delta_E.
  \]

\begin{theorem}\label{thm:OS_basis}
  Let $\M$ be a matroid.  Then
  \[
  \os(\M)=\OLambda_E/
  \bigcap_{B\in\B(\M)}\angl{\Delta_{F(B)}}.
  \]
\end{theorem}
\begin{proof}
  Corollary~3.5 of \cite{DY01} implies that
  $I_{\M}\supseteq \bigcap_B\angl{\Delta_{F(B)}}$.  It is sufficient, then,
  to show that each ideal generator $\partial(e_C)$ is contained in each
  ideal $\angl{\Delta_{F(B)}}$.  Accordingly, if $C$ is a circuit and $B$
  is a basis, let $i$ be the least index for which $C\subseteq F(B)_i$.
  We claim $\abs{C\setminus F(B)_{i-1}}\geq2$: if not, by the choice of
  $i$, we have $C\setminus F(B)_{i-1}=\set{j}$ for some element $j$.
  Since $C$ is a circuit, the span of $C\setminus j$ contains $j$.  Since
  $F(B)_{i-1}$ is a closed set, this implies $j\in F(B)_{i-1}$, a contradiction.

  It follows that $C\setminus F(B)_{i-1}\supseteq\set{j,j'}$ for some
  distinct $j$, $j'$, so $\partial(e_C)$ is a multiple of $e_j-e_{j'}\in
  \Delta_{F(B)}$.  Then $\partial(e_C)\in\angl{\Delta_{F(B)}}$, as required.
 \end{proof}
\subsection{Activities and broken circuits}\label{subsec:BC}

The independent sets of $\M$ form a simplicial
complex $\IN(\M)$ on the set $E$.  If $B$ is a basis for $\M$, an element
$e\not\in B$ is \newword{externally active} for $B$ if $e=\min C(B,e)$, where
$C(B,e)$ is the (unique) circuit in $B\cup\set{e}$.  An element $e\in B$ is
\newword{internally active} for $B$ if $e=\min C^*(B,e)$, where $C^*(B,e)$ is
the (unique) cocircuit in $B$ containing $e$.  We write $\EA(B)$ and
$\IA(B)$ for the set of externally and internally active elements for $B$,
respectively.  For an independent set $I$, we define $\EA(I)$ to be the
externally active elements for $I$ regarded as a basis of $\M\vert_{\cl(I)}$.

The sets $I\in\IN(\M)$ with $\EA(I)=\varnothing$ form a subcomplex of $\IN(\M)$,
the \newword{broken circuit complex}, denoted $\BC(\M)$.
It is a pure simplicial complex of dimension
$r-1$.  The element $e_0:=\min E$ is in every
maximal face, so
$\BC(\M)$ is the cone over the induced subcomplex on vertices $E\setminus
\set{e_0}$, denoted $\bc(\M)$ and called the \newword{reduced broken circuit
complex}.  For details, we refer to
the survey of Bj\"orner~\cite{Bj92}.

Broken circuits are of importance because the elements $\partial(e_C)\in I_{\M}$
indexed by circuits form a universal Gr\"obner basis.  Ordering suitably,
the initial ideal $\In(I_{\M})$
consists of monomials $e_I$ for which $I$ contains a broken circuit.  That is,
if we
let $\k\angl{K}$ denote the exterior Stanley--Reisner ring of a simplicial
complex $K$, 
the Gr\"obner deformation of $\OS(\M)$ is the monomial quotient
\begin{equation}\label{eq:monomial iso}
\Lambda_E/\In(I_{\M})\cong \k\angl{\BC(\M)}.
\end{equation}
Similarly, for $\os(\M)\cong \OS(\M;e_0)$, we obtain the exterior
Stanley--Reisner ring of the reduced broken circuit complex,
$\k\angl{\bc(\M)}$.

The \newword{$\beta$-non-broken-circuits} are defined to be
\begin{equation}\label{eq:bnbc}
\bnbc(\M)=\set{B\in\B(\M)\colon \abs{\IA(B)}=1\text{~and~}\abs{\EA(B)}=0}.
\end{equation}
The name refers to the fact that $\abs{\bnbc(\M)}=\beta(\M)$ (see, e.g.,
\cite[Prop.\ 6.2.12]{BrOx92}.)

\subsection{Resonance}\label{ss:resonance}
The cohomology of Orlik--Solomon algebras in relation to local system
cohomology on arrangement complements has a long history, with origins in
\cite{ESV92,STV95}.  Resonance varieties of arrangements receive
systematic study beginning in \cite{LY00,FY07}.  Resonance varieties of
exterior Stanley--Reisner rings are more easily characterized,
and a complete description is due to Papadima and Suciu~\cite{PS09}, using
earlier work of Aramova, Avramov and Herzog~\cite{AAH00}.

We consider resonance in the following setting.  For a graded module $A$ over an
exterior algebra $\Lambda V$, multiplication by an element $w\in V$
gives $A$ the structure of a cochain complex $(A^\bfcdot,w)$.  For integers $p$ and
$d\geq0$,
\[
\R_d^p(A):=\set{w\in V\colon\dim_{\k} H^p(A^\bfcdot,w)\geq d}.
\]
are the \newword{resonance varieties} of $A$. We abbreviate $\R^p(A):=
\R^p_1(A)$.
We write $\R^p_d(\M):=\R^p_d(\os(\M))$, the (reduced) \newword{resonance
  varieties} of $\M$.  At least for
$d=1$ it is sufficient to consider the reduced
Orlik--Solomon algebra:
\begin{proposition}
  For any $\M$ and $0\leq p\leq \rank(\M)-1$, we have
  $\R^p(\OS(\M))=\R^p(\os(\M))$.
\end{proposition}
\begin{proof}
  By \cite[Prop.\ 3.7]{De16} we have $\R^p(\OS(\M))=\R^p(\os(\M))\cup
  \R^{p-1}(\os(\M))$, and $\R^{p-1}(\os(\M))\subseteq\R^p(\os(\M))$ by
  \cite[Thm.\ 4.1(b)]{EPY03}.  
\end{proof}
We write $\L^{\irr}(\M)$ to denote the collection of connected flats of $\M$,
and $\L^{\irr}_{\leq p}(\M)$ for those of rank at most $p$.
The following is an immediate consequence of the definitions.
\begin{proposition}\label{prop:cone}
  If a matroid $\M$ satisfies the convolution
  resonance formula~\eqref{eq:res formula}, then
  \[
  \R^p(\M)=\bigcup_{F\in\L_{\leq p+1}^\irr(\M)}\Delta_F
  \]
\end{proposition}
The corresponding expression for $\R^p_d(\M)$ is not particularly illuminating;
however, the Hilbert series of the graded module $H^\bfcdot(\os(\M),w)$ can
be expressed in terms of the Tutte polynomial.  For each nonzero $w\in \Delta_E$, if $\M$ satisfies \eqref{eq:res formula}, then 
\[
h(H^\bfcdot(\os(\M),w),t)=\sum_{\substack{
  F\in \L^\irr(\M)\colon\\
  w\in \Delta_F}}
\beta(\M|_F)\pi(\M/F;t)t^{\rank(F)-1},
\]
where $\pi(\M;t):=t^r\cdot T_{\M}(t^{-1}+1,0)$ is the Poincar\'e polynomial.
Reindexing the sums, for nonzero $w$, we have
\begin{align}
\sum_{0\leq p\leq r-1}\dim_\k H^{r-1-p}(\os(\M),w)t^p&=
\sum_{\substack{
  F\in \Z(\M)\colon\\
  w\in \Delta_F}} \beta(\M|_F)T_{\M/F}(t+1,0)\nonumber\\
&=\sum_{\substack{
  F\in \Z(\M)\colon\\
  w\in \Delta_F}} \partial_yT_{\M|_F}(0,y)|_{y=0}T_{\M/F}(t+1,0),\label{eq:conv}
\end{align}
since nonempty cyclic flats contain at least two elements,
and the $\beta$-invariant is the coefficient of $y$
\cite[Prop.\ 6.2.12]{BrOx92}.  We single out a special case:
\begin{proposition}\label{prop:tutte}
  If $\M$ satisfies the convolution resonance formula and
  $w\in \Delta_F$ for all nonempty cyclic flats $F$, then
  \[
  \sum_{0\leq p\leq r-1}\dim_\k H^{r-1-p}(\os(\M),w)t^p=
  \partial_y T_{\M}(t+1,y)|_{y=0}
  \]
\end{proposition}
\begin{proof}
  We apply the convolution formula of Kook--Reiner--Stanton~\cite[Thm.\ 1]{KRS99} to the
  right-hand side of \eqref{eq:conv}.  Their sum is over all subsets; however,
  we note that $T_{\M/S}(t+1,0)=0$ unless $S$ is a flat,
  $\beta(\M|_S)=0$ unless $S$ is connected, and (for a connected flat $S$),
  $w\not\in \Delta_S$ unless
    $S$ is cyclic.
\end{proof}

\section{Special classes of matroids}
\subsection{Arboreal matroids}\label{ss:arboreal}
Now we return to the arboreal matroids defined in the
introduction.
We will call a tree $T$ on $E$ a \emph{flat tree} if
it restricts to a subtree on each element of $\L^\irr_{\geq1}(\M)$.
By Lemma~\ref{lem:partial_product}, $T$ has the property that, for any flat
$F\in\L^\irr_{\geq1}(\M)$, the set $\bB(T[F])$ is a basis for $\Delta_F$,
where $T[F]$ denotes the induced subgraph of $T$ on $F$.

\begin{lemma}\label{lem:interval}
  Let $\M$ be a matroid on $E$. If there exists a total order on $E$ for
  which every connected flat of $\M$ is an interval, then $\M$ is arboreal.
\end{lemma}
\begin{proof}
The path on $E$ induced by such an order is a flat tree $T$ of $\M$.
\end{proof}
\begin{example}\label{ex:arboreal}
  All lattice path matroids are arboreal, since their connected flats
  are intervals, by \cite[Thm.\ 5.7]{BdM06}.  The class of laminar matroids
  of Fife and Oxley~\cite{fifeOxley} is arboreal (\cite{U26+}).  All rank-$2$
  matroids are arboreal.

  On the other hand, the $3$-whirl, the underlying
  matroid of the arrangement defined by the polynomial
  $Q=xyz(x+y)(x+z)(y+z)$, is not arboreal: it has
  rank-$2$ cyclic flats $013$, $024$, and $125$.  The cyclic flats have
  nonempty pairwise intersections, but $013\cap 024\cap 125=\emptyset$,
  so $\L^\irr_{\geq1}(\M)$ does not satisfy the Helly property.
\end{example}
Now we prove that, if $\M$ is arboreal, its Orlik--Solomon algebra is an
exterior face ring.  For this, we need to examine the order of the ground
set.  If $T$ is a tree, the (vertex) sets of subtrees forms a convex
geometry, in the sense of Edelman and Jamison~\cite{EJ85}, who define a
compatible ordering $e_0<e_1<\cdots<e_n$ of $E$ to be one in which
the span of each subset $\set{e_0,\ldots,e_i}$ for $0\leq i\leq n$ is a
subtree.  Clearly, compatible orderings exist for all trees $T$.  In
particular, if $T$ is a path, its natural order is compatible, since
initial segments are subpaths.
With this detail in mind, we can state Theorem~\ref{thm:degeneration} more
precisely:
\begin{theorem}\label{thm:arboreal_OS}
  If $\M$ is an arboreal matroid with flat tree $T$ and $E$ is compatibly
  ordered with respect to the convex geometry of subtrees, then
  \[
  \OS(\M)\cong\k\angl{\BC(\M)}.
  \]
\end{theorem}
\begin{proof}
  By relabelling, we assume the standard order on $E=[n]$ is compatible.
  It is enough to show that $\os(\M)\cong\k\angl{\bc(\M)}$.
  
  The vertices $[i]$ span a subtree of $T$ for which $i$ is a leaf, for each $1\leq i\leq n$. Let $\ell_i$ denote the (unique)
  vertex in $[i-1]$ with an edge to 
  $i$. We remark that $\ell_i\in F$ for $i\in F\setminus \min F$ since $F\cap [i]$ is a subtree on at least two vertices with leaf $i$. Let $b_i:=e_i-e_{\ell_i}$.  If we direct the
  edges of $T$ using the order on $E$, the vectors
  $\bB(T)=\set{b_1,\ldots,b_n}$ constitute a basis of $\Delta_E$
  (Lemma~\ref{lem:partial_product}).

  If $F$ is a (nonempty) flat of $\M$, tree-convexity implies that the
  induced subgraph $T[F]$ is also a tree, which means $\bB(T)$
  restricts to a basis of $\Delta_F$,
  \[
  \bB(T[F]) = \set{b_i\colon i\in F\setminus \min(F)}.
  \]
  Theorem~\ref{thm:OS_pres} now presents $\os(\M)$ as a quotient by a monomial
  ideal in the variables $\bB(T)$,
  \[
  I_{\M}=\sum_{F\in\L^\irr_{\geq1}(\M)}\angl{b_S\colon S\subseteq F\setminus\min(F),\,
    \abs{S}=\rank(F)}.
  \]

  It remains to show that $I_{\M}$ is the ideal of non-faces of $\bc(\M)$,
  which we denote $I_{\bc}$.
  For this,
  choose a generator $b_S$ of $I_{\M}$.  Since $\abs{S\cup\min(F)}>\rank(F)$,
  there is a circuit $C\subseteq S\cup \min(F)$.  If $C\subseteq S$, then
  $S$ contains a broken circuit $C\setminus\min(C)$, and $b_S\in I_{\bc}$.
  If $C\not\subseteq S$, then $\min C=\min F$, and again $S$ contains a
  broken circuit, $C\setminus\min(F)$; either way, $b_S\in I_{\bc}$.

  In the other direction, choose a generator $b_S\in I_{\bc}$.  Then
  $S=C\setminus\min(C)$ for some circuit $C$: let $F=\cl(C)$.  Then
  $\abs{S}=\rank(C)=\rank(F)$.
  Since
  $\min(F)\leq\min(C)$, we also have $S\subseteq F\setminus\min(F)$.  These
  facts taken together show $b_S$ is a generator of $I_{\M}$.
\end{proof}

\begin{remark}
  It should be noted that the isomorphism above does not hold for
  arbitrary matroids.  For example, for the matroid $\M$ of the complete graph
  $K_4$, one can distinguish $\OS(\M)$ from $\k\angl{\BC(\M)}$ using
  resonance.  It is well-known that $\R^1(\M)$ consists of $5$ two-dimensional
  components (e.g., \cite[Ex.\ 10.3]{Su01}), while $\R^1(\k\angl{\BC(\M)})$
  has only four (using \cite[Prop.\ 3.6]{PS09}.)
\end{remark}

\subsection{Schubert matroids}
Let $E=\set{\varepsilon_1<\varepsilon_2<\ldots<\varepsilon_n}$ be an
ordered ground set, and $S\subseteq E$ a subset of size $r$.  The
\newword{Schubert matroid} $\M_{S,E}$ is a matroid on $E$ of rank $r$ with
several equivalent formulations.  
(We write $\M_S$ when the choice of $E$ is clear.)  We briefly recall some of
these.

The \newword{Gale order} is defined by letting $S\preceq T$ if and
only if $S=\set{s_1<\cdots<s_r}$ and $T=\set{t_1<\cdots<t_r}$ satisfy
$s_i\leq t_i$ for all $1\leq i\leq r$.
The \newword{Schubert matroid} $\M_S$ has bases
\[
\B(\M_S):=\set{B\subseteq E\colon \abs{B}=r\text{~and~}B\succeq S}.
\]
From this point of view, $\M$ is a Schubert matroid if and only if
$\IN(\M)$ is a pure, shifted complex: see~\cite{Kl03},~\cite[Cor.\ 4.4]{Ard03}.

\begin{proposition}[Lem.\ 2,~\cite{OPR82}; see also Lem.\ 2.2,~\cite{BdM08}]
  \label{prop:Schubert_cyclic}
  Suppose $\abs{E}\geq1$ is ordered.  Then $\M$ is a Schubert matroid
  if and
  only if its cyclic flats $\Z(\M)$ are all of the form $E_{\leq z}$,
  for some $z\in E$.

  If $\M=\M_S$ is a Schubert matroid, each cyclic flat $E_{\leq z}$
  has rank $\abs{S_{\leq z}}$.\label{eq:fact1}  In this case,
  let $Z(S,E):=\set{z\in E\colon E_{\leq z}\in\Z(\M)}$.
\end{proposition}

Schubert matroids are also obtained by iteratively adding coloops and
performing free extensions.  This was
Crapo's original formulation~\cite{Cr65}: see the discussion just prior
to~\cite[Exm.\ 6.6]{crapo-schmitt}.
\begin{proposition}\label{prop:Schubert_facts2}
  Suppose $\abs{E}\geq 1$ and $S\subseteq E$.
\begin{enumerate}
    \item
      Let $\varepsilon=\max E$, then let $S'=S\setminus\varepsilon$ and
      $E'=E\setminus\varepsilon$.
      If $\varepsilon\in S$, then $\M_{S,E}=\M_{S',E'}\oplus \U_{1,1}$.  If
      $\varepsilon\not\in S$, then $\M_{S,E}=\M_{S',E'}+\varepsilon$, a free extension. \label{eq:facts2(1)}
    \item 
      Let $\varepsilon=\min E$, then let $S'=S\setminus\varepsilon$ and
      $E'=E\setminus\varepsilon$.
      If $\varepsilon\not\in S$, then $\M_{S,E}=\U_{0,1}\oplus \M_{S',E'}$.
      If $\varepsilon\in S$,
      then $\M_{S,E}=\varepsilon\times\M_{S',E'}$, a free coextension.\label{eq:fact3'}
  \end{enumerate}
\end{proposition}

We regard $\M_S$ as a \newword{lattice path matroid} as follows
(for details, see \cite{BdMN03,BdM06}).  A subset $T\subseteq E$ indexes
a \newword{lattice path} of length $n$ in $\RR^2$ from $(0,0)$ to
$(n-\abs{T},\abs{T})$ where
the $i$th step goes up if $\varepsilon_i\in T$ and the $i$th step goes right
otherwise.  We will freely make use of this identification of subsets and
paths from now on.  If a path has $r$ vertical steps, those steps
form a basis if and only if
the path lies (weakly) to the right of the path of $S$.

\begin{figure}
\input{sch_arxiv.fig1.tex}
\caption{$I=245$ is independent in $\M_{1346,[1,7]}$,
  but $3\cup I$ is not, since its path exits $\lambda_{S,E}$.
  \label{fig:LPM}}
\end{figure}
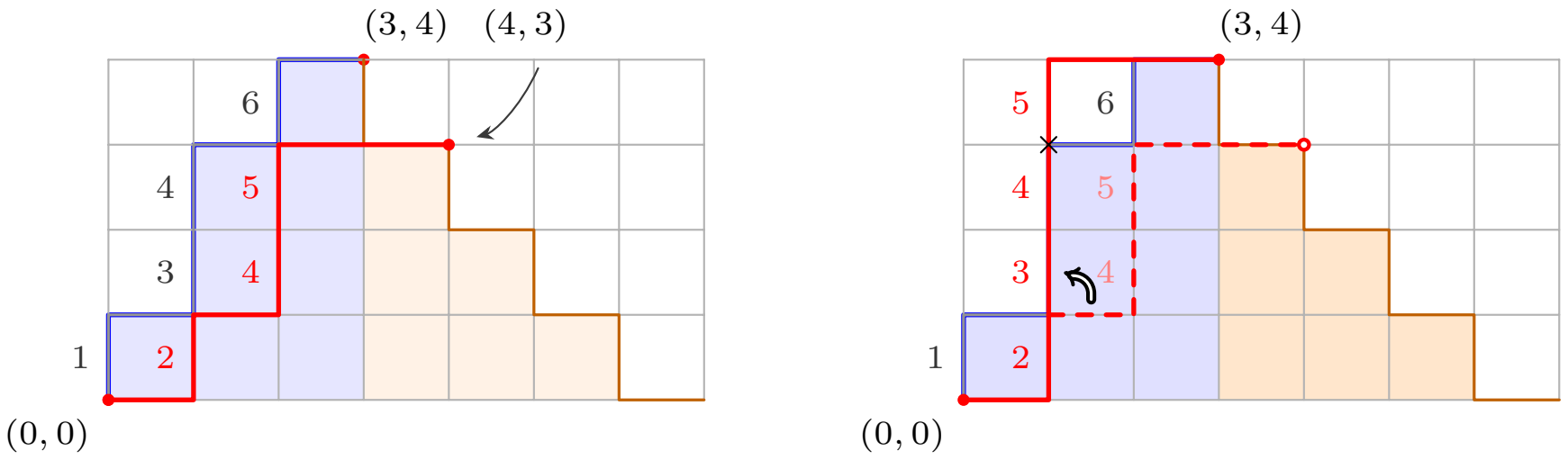

More generally, $I\subseteq E$ is independent if and only if its lattice
path lies in the partition shape determined by the first $\abs{I}$-many
elements of $S$.  The independent sets can then
be characterized as those lattice paths of length $n$ from $(0,0)$ that
lie inside a figure $\lambda_{S,E}$ bounded
by the path $S$, a staircase through vertices $(n-k,k)$ for $0\leq k\leq r$,
and the $y$-axis.  See Figure~\ref{fig:LPM}.

We can recognize the closure of an independent set in terms of lattice paths:
\begin{propdef}\label{def:tight}
  Suppose $I$ is an independent set in $\M_S$.  Say
  $\varepsilon\in E\setminus I$ is \newword{tight} if $\varepsilon$ is in the
  closure of $I$; otherwise, say $\varepsilon$ is \newword{slack}.
  An element $\varepsilon$ is slack if and only if the lattice path of
  $I\cup\varepsilon$ is contained in $\lambda_{S,E}$.
\end{propdef}
\begin{proof}
  We have $\varepsilon\not\in\cl(I)$ if and only if $I\cup\varepsilon$
  is independent, if and only if its path remains in $\lambda_{S,E}$.
\end{proof}
\begin{example}\label{ex:indep_set}
  For the set $I=245$ shown in Figure~\ref{fig:LPM}, steps $1$ and $3$
  are tight, while steps $6$ and $7$ are slack.
\end{example}

\begin{proposition}\label{prop:tightinitial}
  For $I\in\IN(\M_{S,E})$, whenever $\varepsilon_j,\varepsilon_{j+1}\not\in I$, if $\varepsilon_{j+1}$ is tight,
  then so is $\varepsilon_j$.
\end{proposition}
\begin{proof}
  Write $I\cup\set{\varepsilon_{j+1}}=
  \set{t_1<\cdots<t_{k+1}}$ in increasing order.
  By hypothesis, $I\cup\set{\varepsilon_{j+1}}\not\succeq S':=\set{s_1<\cdots<s_{k+1}}$, so pick some $i$ for which $s_i>t_i$. There are two cases to consider: if $t_i=\varepsilon_{j+1}$, then since $\varepsilon_{j+1}\not\in I$
  we have $s_i>\varepsilon_j$ as well, so $I\cup\set{\varepsilon_{j}}
  \not \succeq S'$; that is, $\varepsilon_{j}$ is tight.
  Otherwise, if $t_i\neq \varepsilon_{j+1}$, then $t_i$ remains
  the $i$th smallest element
  of $I\cup\set{\varepsilon_{j}}$, and again $I\cup\set{\varepsilon_{j}}
  \not\succeq S'$.
\end{proof}

Cyclic flats can be understood in terms of lattice paths as well.
\begin{propdef}\label{def:corner}
  Given a set $E$, introduce elements $\varepsilon_0$ and $\varepsilon_{n+1}$ with
  $\varepsilon_0<\min E$ and $\max E<\varepsilon_{n+1}$.
  We say $\varepsilon_i\in \set{\varepsilon_0}\cup E$ is an
  \newword{inside corner} of $\M_{S,E}$ if
  $\varepsilon_i\not\in S$ and
  $\varepsilon_{i+1}\in S\cup\set{\varepsilon_{n+1}}$.  The set of inside
  corners of $\M_{S,E}$ equals the set $Z(S,E)$ of maximal elements of cyclic
  flats.
\end{propdef}

\begin{remark}\label{rem:extremes}
  We note that $\varepsilon_0$ is an inside corner if and only if $\min S=
\min E$: i.e., its path starts with a vertical step.  The last element of $E$
is an inside corner if and only if $\max S<\max E$: i.e., it ends with a
horizontal step.
\end{remark}

\begin{example}\label{ex:1346}
  Let $S=1346\subseteq E=[1,7]$.
  Then $\M_S$ has cyclic flats
  $\varnothing\subset [1,2] \subset [1,5] \subset E$. Noting $\min S = \min E$, the inside corners
  $Z(S,E)=\set{0,2,5,7}$ 
  are in positions $(0,0)$, $(1,1)$, $(2,3)$, and $(3,4)$;
  see Figure~\ref{fig:LPM}. 
\end{example}

\begin{notation}\label{not:rS}
  When $\varepsilon=\min S=\min E$, Proposition~\ref{prop:Schubert_facts2}\eqref{eq:fact3'}
  gives $\M_{S,E} = \varepsilon\times \M_{\rS,\rE}$, where $\rS:=S\setminus\varepsilon$ and $\rE:=E\setminus\varepsilon$.
We see that $\lambda_{\rS,\rE}$ is obtained from $\lambda_{S,E}$ by removing
the bottom row.
\end{notation}

We end this section by recording a useful fact about restriction.
\begin{lemma}[Thm.\ 3.1, \cite{BdM06}]\label{lem:Schubert_restr}
  Suppose $P=E_{\leq \varepsilon}$ for some $\varepsilon\in E$.  Then
  \[
  \M_{S,E}\vert_P = \M_{S_{\leq \varepsilon},P}.
  \]
\end{lemma}

\subsection{Orlik--Solomon algebras of Schubert matroids}
We noted in the previous section that Schubert matroids are
lattice path matroids, hence they are arboreal (Example~\ref{ex:arboreal}).
By Theorem~\ref{thm:arboreal_OS}, their Orlik--Solomon algebras are face
rings.  Our assumption that $\M_{S,E}$ is loopless is equivalent to the
condition that $\min S = \min E$.

\begin{proposition}\label{prop:OS_IN}
  For a Schubert matroid, $\os(\M_S)\cong \k
  \angl{\IN(\M_{\rS,\rE})}$.
\end{proposition}
\begin{proof}
  Brylawski~\cite[Thm.\ 4.2]{Bry77} showed that, if $\M$ is a matroid on
  $E$ and $\varepsilon<\min E$, then $\bc(\varepsilon\times \M)$
  equals $\IN(\M)$.  The claim then follows from
  Proposition~\ref{prop:Schubert_facts2}\eqref{eq:fact3'} and
  Theorem~\ref{thm:arboreal_OS}.
\end{proof}

\subsection{Activities for Schubert matroids}
Basis activities have simple lattice path descriptions.
\begin{theorem}[Thm.\ 5.4, \cite{BdMN03}]\label{thm:Bactive}
  For a basis $B\in \B(\M_{S,E})$, the set $\EA(B)$ consists of horizontal steps
  in the path of $B$ on the bottom edge of $\lambda_{S,E}$, while
  $\IA(B)$ consists of vertical steps in the path of $B$ on the left edge
  of $\lambda_{S,E}$.
\end{theorem}
This gives a lattice-path interpretation of Proposition~\ref{prop:OS_IN}:
$\nbc$ bases have no externally active elements, so they are paths that
start with a vertical step.  Paths in $\lambda_{S,E}$
that start with a vertical step are in bijection with the set of
all paths in $\lambda_{\rS,\rE}$.

\begin{corollary}\label{cor:bnbc}
  Let $\varepsilon=\min E$.  The map $B\mapsto B\setminus\varepsilon$
  is a bijection between $\bnbc(\M_{S,E})$ and
  \[
  \set{B\in\B(\M_{\rS,\rE})\colon \IA(B)=\varnothing}.
  \]
\end{corollary}
\begin{proof}
  A basis $B\in\bnbc(\M_{S,E})$ has $\IA(B)=
  \set{\varepsilon}$, since as a lattice path, only its first step $\varepsilon$
  lies on the left edge of $\lambda_{S,E}$.  Removing the bottom row of
  $\lambda_{S,E}$, then, leaves a lattice path in $\lambda_{\rS,\rE}$
  that avoids the left edge, and this process is easily seen to be reversible.
\end{proof}

We need to extend Theorem~\ref{thm:Bactive} to submaximal independent
sets.  
\begin{proposition}\label{prop:Iactive}
  For any $I\in\IN(\M_{S,E})$, the set $\EA(I)$ consists of all tight steps
  $\varepsilon$ in the path of $I$ for which $\varepsilon<\min I$.
\end{proposition}

\begin{proof}
 Write $\M:=\M_{S,E}$ and let $F:=\cl_\M(I)$,
  the union of $I$ and tight steps $\varepsilon$. As $F$ contains a cyclic flat since it contains a circuit of $\M$, we note $\Z(\M\mid_F)\subseteq \Z(\M)$ as a subchain; thus, the restriction $\M\mid_F$ is again a Schubert matroid, say with some Gale support $S_F$. As $I$ is a basis of $\M\mid_F$, we note $I\succeq S_F$. By Theorem~\ref{thm:Bactive}, externally active elements of basis $I$ of $\M\mid_F$ are then $\varepsilon\in F$ for which $\varepsilon<\min(I)$.
\end{proof}

To see that this recovers Theorem~\ref{thm:Bactive} when $I$ is maximal,
we observe that all horizontal steps $\varepsilon$ in a basis are
tight, so the condition that $\varepsilon<\min I$ is the same as
$\varepsilon$ lying on the bottom edge of $\lambda_{S,E}$.

\begin{example}[Example~\ref{ex:indep_set}, continued]
    For $I=245$ shown in Figure~\ref{fig:LPM},
    $\cl(I)=[1,5]$, and $\EA(I)=1$.
\end{example}

\section{Resonance for Schubert arrangements}
Our goal in this section is to prove Theorem~\ref{thm:schres}.  We start
by noticing the Orlik--Solomon algebra of a Schubert matroid has a large
automorphism group.  In this section, we assume $E=[n]$ for some
integer $n\geq1$.

Given a set of integers $Z=\set{z_1<\cdots<z_k}$, the flag of subspaces
\begin{equation}\label{eq:Zchain}
\Delta_{[z_1]}\subseteq \Delta_{[z_2]}\subseteq \cdots \subseteq
\Delta_{[z_k]}\subseteq \Delta_E
\end{equation}
is stabilized by a parabolic subgroup of the general linear group,
\[
P_Z := \set{g\in \GL(\Delta_E)\colon g(\Delta_{[z_i]})\subseteq
\Delta_{[z_i]}\text{~for all $1\leq i\leq k$}}.
\]
\begin{lemma}\label{lem:auts}
  The automorphism group of $\os(\M_S)$ contains $P_{Z(S,E)}$.
\end{lemma}
\begin{proof}
  The set $\set{[z]\colon z\in Z(S,E)}$ is the set of cyclic flats
  (Proposition~\ref{prop:Schubert_cyclic}), so it contains all 
  connected flats of $\M_S$.  It follows that
  $P_{Z(S,E)}$ stabilizes the defining ideal of $\os(\M_S)$, using the
  description provided by Theorem~\ref{thm:OS_pres}.
\end{proof}
\subsection{Reduction to monomials}
These automorphisms allow us to reduce the calculation of resonance to that
of an exterior
Stanley--Reisner algebra.  For each $w\in\k^E$, we let
\begin{equation}\label{eq:zw}
  z(S,E;w)=\min(\set{z\in Z(S,E)\colon w\in\Delta_{[z]}}\cup\set\infty).
\end{equation}
\begin{theorem}\label{thm:reduction}
  For each nonzero $w\in\Delta_E$, if $z:=z(S,E;w)<\infty$, we have
  \[
  H^\bfcdot(\os(\M_S),w)\cong H^\bfcdot(\k\angl{\IN(\M_{\rS,\rE})},e_z).
  \]
  If $z=\infty$, then $H^\bfcdot(\os(\M_S),w)=0$.
\end{theorem}
\begin{proof}
  First, suppose $z:=z(S,E;w)<\infty$, so that 
  $w\in\Delta_{[z]}\setminus \Delta_{[z']}$, where $z'$ is the 
  preceding inside corner.
  Since $w':=e_{z}-e_{z-1}\in\Delta_{[z]}\setminus \Delta_{[z']}$
  as well, we have $w'=g\cdot w$ for some $g\in P_{Z(S,E)}$. Since a connected flat is always cyclic, we apply Lemma~\ref{lem:auts}, noting that
  the standard order on $E$ gives a flat tree.  The result now follows by applying Proposition~\ref{prop:OS_IN}.

  Otherwise, $z=\infty$: since $w\in\Delta_E$, the set
  $E$ is not a cyclic flat. That is, $n':=\max S<\max E$, and
  $\M_{S,E}=\M_{S,E'}\oplus \U_{k,k}$ where $k=n-n'$ and $E'=E_{\leq n'}$,
  using Proposition~\ref{prop:Schubert_facts2}\eqref{eq:facts2(1)}.
  Orlik--Solomon algebras
  of sums are tensor products, giving an isomorphism of cochain complexes
  \[
  (\os(\M_{S,E}),w)\cong 
  (\OS(\M_{S,E};0),\phi_0^{-1}(w))\cong(\OS(\M_{S,E'};0),w_1)\otimes_\k(\OS(\U_{k,k}),w_2),
  \]
  where $\phi_0$ is the isomorphism from Proposition~\ref{prop:loops}, and 
  $\phi_0^{-1}(w)= w_1+w_2$ for some $w_1$, $w_2$.  Since
  $z> n'$, we have $w\not\in\Delta_{E'}$, which means $w_2\neq 0$.
  Since $\OS(\U_{k,k})$ is an exterior algebra, the second factor is exact
  when $w_2\neq 0$, and that implies $(\os(\M_{S,E}),w)$ is exact as well.

\end{proof}
\begin{example}[Example~\ref{ex:indep_set}, continued]
  We saw the inside corners $Z(S,E)$ of $\M_{1346,[1,7]}$
  consist of $\set{0,2,5,7}$ (Figure~\ref{fig:LPM}).
  The element $w:=e_1+e_2-e_3-e_4\in
  \Delta_{[1,5]}$ but $w\not\in\Delta_{[1,2]}=\k\text{-}\mathrm{span}\set{e_2-e_1}$, so $w=g\cdot(e_5-e_4)$ for
  some automorphism $g$, and
  \[
  H^\bfcdot(\os(\M_S),w)\cong H^\bfcdot(\os(\M_S),e_5-e_4)=
  H^\bfcdot(\k\angl{\IN(\M_{346,[2,7]})},e_5).
  \]
\end{example}

\begin{proposition}\label{prop:link}
  Suppose $z$ is an inside corner for $\M_{S,E}$, 
  and let $K:=\IN(\M_{\rS,\rE})$.
  Then
  a basis for $H^\bfcdot(\k\angl{K},e_z)$ is given by monomials
  \[
  \set{e_I\colon z\not\in I\in K\setminus\link_K({z})}.
  \]
\end{proposition}
\begin{proof}
  Since $\k\angl{K}$ is a quotient of an exterior algebra by a monomial ideal, it is evident that
  $H^\bfcdot(\k\angl{K},e_z)=\ann(e_z)/(e_z)$
  has a monomial basis consisting of some monomials $e_I$ where
  $I\in K$ and $z\not\in I$.  We observe $e_I\in\ann(e_z)$ exactly when
  $I\cup\set{z}\not\in K$, which is to say $I\not\in\link_K(\set{z})$.
\end{proof}

The monomial basis above has a straightforward lattice path description.
The shape $\overline{\lambda}:=
\lambda_{\rS,\rE}$ is obtained by removing the bottom row.
A lattice path $I$ in $K$ with $z\not\in I$
is one that is
contained in $\overline{\lambda}$ and travels right at step $z$.  The
condition that $I\not\in\link_K(z)$ states that $I\cup z$ is dependent: that
is, $I$ indexes a monomial above if and only if $z$ is a tight edge in $I$
(Proposition~\ref{def:tight}).

\begin{figure}
\input{sch_arxiv.fig2.tex}
\caption{$I=1346$ and $\sh_5(I)=2456$
  in $\M_{1346,[1,7]}$\label{fig:shift}}
\end{figure}
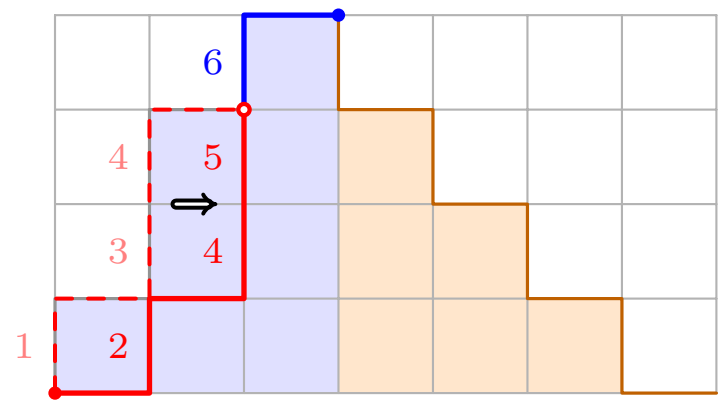

Let $\sh_z\colon E\to E$ be the \newword{shift map} defined by
\[
\sh_z(j) = \begin{cases} j+1 & \text{if $j<z$};\\
  j& \text{otherwise.}
\end{cases}
\]
The shift map extends to a bijection $\sh_z\colon
2^{E\setminus z}\to 2^{E\setminus 0}$.  Whenever 
$z\not\in I$, the lattice path of $\sh_z(I)$ is obtained from the
lattice path of $I$ by adding a horizontal step at the start, and removing
the horizontal step $z$: see Figure~\ref{fig:shift}.

\begin{figure}
  \input{sch_arxiv.fig3.tex}
\caption{Paths meeting an inside corner\label{fig:cases}}
\end{figure}
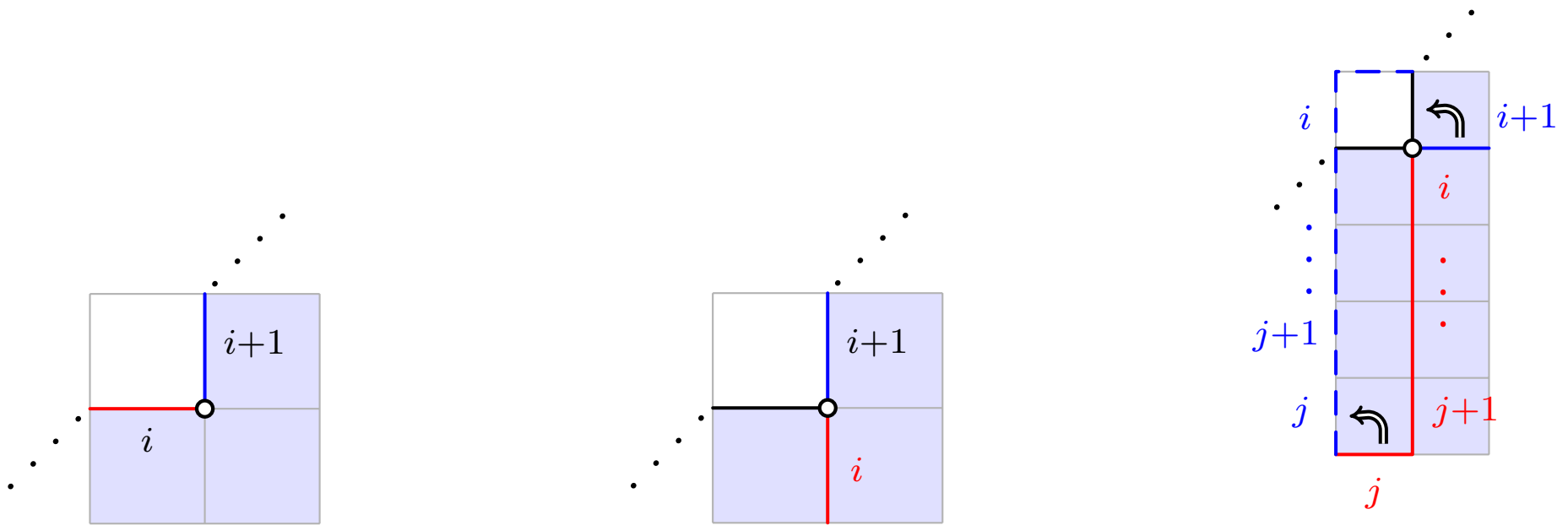

\begin{lemma}\label{lem:key}
  Consider a path $I$ in a partition diagram $\lambda:=\lambda_{S,E}$.
  Let $i$ be the largest element of $E$ for which each vertical step of
  $I_{\leq i}$ lies inside the interior of $\lambda$, and for which
  each horizontal step
  of $I_{\leq i}$ is tight.  Then $i$ is an inside corner.
\end{lemma}
\begin{proof}
  If $i<\max E$, we consider step $i+1$. Suppose first it is vertical:
  by maximality, it lies on the left edge.  If (the previous) 
  step $i$ is horizontal, we obtain an inside corner.
  (Figure~\ref{fig:case1}.)
  If step $i$ is vertical, it lies in the interior of $\lambda$,
  so again it gives an inside corner (Figure~\ref{fig:case2}.)

  Now suppose $i+1$ is horizontal, hence slack.  First we consider the case
  where $i$ is horizontal.  Slack edges lie in the interior, so we show that
  edge $i$ is not in the interior in order to establish there is an inside
  corner at step $i$.

  For this, order the elements
  $I\cup\set{i+1}:=\set{t_1<\dots<t_{k+1}}$ and suppose $t_p=i+1$.  Let $S':=\set{s_1<\dots<s_{k+1}}$, the first $k+1$ elements of $S$.
  Then $I\cup\set{i+1}\succeq S'$. 

  Because $i$ is horizontal, it must be tight, which implies
  $i<s_p\leq t_p=i+1$: that is, $s_p=i+1\in S'$.
  We claim $i\not\in S'$: if it were, we would have $i=s_{p-1}\leq t_{p-1}
  <t_p=i+1$, forcing $t_{p-1}=i$, hence $i\in I$, which would make $i$ a
  vertical step instead. Since $i+1\in S$, the step $i$ is an inside corner. 

  Otherwise, $i$ is vertical (hence in the interior.)
  We claim its endpoint must lie on the left edge.
  For this, suppose $j<i$ is horizontal and $j+1,j+2,\ldots,i$ are
  vertical (interior edges).  By maximality, $j$ is tight, so $I\cup j\not
  \succeq S'$, while $I\cup (i+1)\succeq S'$.  Therefore 
  the path $I\cup j$ leaves $\lambda$ at step $i$, which is to say that
  $i$ is an inside corner (Figure~\ref{fig:case3}.)

  If $i=\max E$, we claim the last step in the path of $S$ is horizontal.
  Again, we check two cases.  If $i\in I$ (a vertical step), then $i\not\in S$,
  since otherwise $i$ would not be in the interior.
  If $i\not\in I$ (a horizontal step), we must also have $i\not\in S$, else
  $I\not\succeq S$.  Either way, $i$ is an inside corner, in view of
  Remark~\ref{rem:extremes}.  
\end{proof}

This leads to the combinatorial result underpinning Theorem~\ref{thm:schres}.
We extend the join of simplicial complexes to arbitrary collections of sets:
if $K,K'\subseteq 2^E$, let
\[
K\star K' := \set{A\cup A'\colon A\in K\text{~and~}A'\in K'}.
\]

\begin{proposition}\label{prop:main_bijection}
  For each inside corner $z$, the map $\Sh_z\colon 2^{\rE\setminus z}\to 2^E$ given by
  $\Sh_z(I)=0\cup\sh_z(I)$ restricts to a bijection
  \begin{equation}\label{eq:pathbijection}
  \set{I\colon z\not\in I\in K\setminus\link_K(z)}\to
  \bigsqcup_{j\in Z(S,E)_{\geq z}}\bnbc(\M\vert_{[j]})\star \BC(\M/ [j])
  \end{equation}
  where $K=\IN(\M_{\rS,\rE})$ and $\M:=\M_{S,E}$.
\end{proposition}
\begin{proof}
  Consider an element $I$ from the left side of
  \eqref{eq:pathbijection}.  Write $J:=\sh_z(I)$ and $m:=\max(J)$. Clearly, $J\succeq
  I\succeq \rS_{\leq m}$,
  so $J$ indexes a lattice path in
  $\overline{\lambda}:=\lambda_{\rS,\rE}$. In fact, $J_{\leq z} \succneqq \rS_{\leq z}$: i.e., each vertical step of $J_{\leq z}$ lies in the interior of $\overline{\lambda}$. As $J_{\leq z}$ is a basis of $\M_{\rS_{\leq z},\rE_{\leq z}}$, each horizontal step is tight. 
       
    Let $j$ be the greatest
  element of $\rE$ so that each vertical
  step of the path $J_{\leq j}$ lies inside the interior of
  $\overline{\lambda}$, and each horizontal step is tight. By construction, $z\leq j$.

  By Lemma~\ref{lem:key}, step $j$ is an inside corner. We consider the partition $J=J_{\leq j}\sqcup J_{>j}$.  By construction,
  $J_{\leq j}$ is a basis in $\M_{\rS_{\leq j},\rE_{\leq j}}$ with no
  internally active elements.  By Corollary~\ref{cor:bnbc}, $0\cup J_{\leq j}$
  is an element of $\bnbc(\M_{S_{\leq j},E_{\leq j}})$.
  Similarly, the set $J_{>j}$
  indexes a lattice path in $\lambda_{S_{> j},E_{> j}}$.  By our choice of
  $j$, step $j+1$ is slack if horizontal.  Using Propositions~\ref{prop:Iactive}
  and \ref{prop:tightinitial}, we conclude that $J_{>j}$ has no
  externally active elements in $M_{S_{>j},E_{>j}}$; that is,
  $J_{>j}\in\BC(\M_{S_{>j},E_{>j}})$.
  Finally, $\M_{S_{\leq j},E_{\leq j}}=\M\vert_{[j]}$ and
  $\M_{S_{> j},E_{> j}}=\M/[j]$.

  To see that $\Sh_z$ is a bijection, we check that the operation
  above is reversible. For an inside corner $j\in Z(S,E)_{\geq z}$,
  consider $J\in \bnbc(\M\vert_{[j]})\star \BC(\M/[j])$. The $\nbc$
  portion satisfies $J_{>j}\succeq \rS\cap (j,m]$, whereas the
    $\bnbc$ portion $J_{\leq j}$ in $\M\vert_{[j]}$ maps uniquely to
    the lattice path $\overline{J}_{\leq j}:=J_{\leq j}\setminus
    0$, which is a basis in $\M_{\rS_{\leq j}, \rE_{\leq j}}$, by the
    bijection of Corollary~\ref{cor:bnbc}. $\overline{J}_{\leq j}$
    lives in the interior of $\lambda_{\rS_{\leq j},\rE_{\leq j}}$ and
    avoids the left edge; thus, $\overline{J}_{\leq
      j}\succneqq \rS_{\leq j}$. Applying the map $\sh^{-1}_{z}:\rE\to E$
\[
i \mapsto \begin{cases} i-1 & \text{if $0<i\leq z$};\\ i&
  \text{otherwise,}
\end{cases}
\]   
to $\bar{J}_{\leq j}\cup J_{>j}$, we observe $\sh^{-1}_z(\bar{J}_{\leq j}\cup
J_{>j})\succeq \rS_{\leq m}$ with tight step $z$. We thus conclude
$J = 0\cup\sh_z(I)$ for some $I\in \IN(\M_{\rS,\rE})$ with tight $z$.
\end{proof}
\begin{example}[Example~\ref{ex:indep_set}, continued]
  We continue the example above with $(\rS,\rE)=(1346,[1,7])$.
  Write $Z:=[5]$. We have $I:=1346\not\in\link(5)$, since $5$
  is tight for $I$.  Each step in the path $\sh_5(I)$ up to $j=5$
  is either an interior vertical edge
  or a tight horizontal one, while step $6$ is neither.  Then the bijection
  in Proposition~\ref{prop:main_bijection} sends $I$ to
  \begin{align*}
    (245,6)& \in\B(\overline{\M}\vert_{\overline{Z}})\times\BC(\overline{\M}/\overline{Z}) \text{~to}\\
    (0245,6)&\in \bnbc(\M\vert_{Z})\times \BC(\M/Z),
  \end{align*}
  noting that the set $245$ has no internally active elements, 
  and $6$ has no externally active elements.
\end{example}
\begin{proof}[Proof of Theorem~\ref{thm:schres}]
  Let $\M=\M_{S,E}$ be a Schubert matroid.  The cyclic flats of positive rank $\Z_{\geq 1}(\M)$ are indexed by inside corners $\varepsilon>0$.
  For any nonzero $w\in \Delta_E$,
  let $z=z(S,E;w)$ from \eqref{eq:zw}.  If $z=\infty$, both sides of
  \eqref{eq:res formula} are zero.  Otherwise,
  \[
  H^\bfcdot(\os(\M),w)\cong H^\bfcdot(\k\angl{\IN(\M_{\rS,\rE})},e_z)
  \]
  by Theorem~\ref{thm:reduction}.  A (graded) basis for the latter is
  provided by monomials from Proposition~\ref{prop:link}.  By
  Proposition~\ref{prop:main_bijection}, these are in bijection with the set
  \[
  \bigsqcup_{j\in Z(S,E)_{\geq z}}\bnbc(\M\vert_{[j]})\star \BC(\M/ [j]).
  \]
  We note $\abs{\bnbc(\M\vert_{[j]})}=\beta(\M\vert_{[j]})$ and sets $0\cup B\in
  \bnbc(\M\vert_{[j]})$ have $\abs{B}=\rank(\M\vert_{[j]})-1$.  The result
  follows by noting that $\nbc$ monomials form an additive basis for
  $\OS(\M/[j])$.
  
  In fact, unravelling the isomorphism from Theorem~\ref{thm:reduction} shows
  that a set $B\sqcup I\in \bnbc(\M\vert_{[j]})\star \BC(\M/ [j])$ maps to a
  cocycle
  \begin{align*}
    \partial(e_B)\cdot e_I&\in H^{\rank(\M\vert[j])-1+\abs{I}}(\os(\M),w')\\
    &\cong H^{\rank(\M\vert[j])-1+\abs{I}}(\os(\M),w),
  \end{align*}
  where $w'=e_z-e_{z-1}$.
\end{proof}
\begin{corollary}
  If $\M_S$ is a Schubert matroid and $w\in\Delta_F$ for all nonempty
  cyclic $F$, then
  \[
  h(H^\bfcdot(\os(\M_S),w),t) = t^{r-1}\partial_y T_{\M_S}(t^{-1}+1,y)|_{y=0}.
  \]
\end{corollary}
\begin{proof}
  Proposition~\ref{prop:tutte} applies to matroids $\M$ that satisfy
  \eqref{eq:res formula}, and its statement is non-vacuous if and only if
  $0\neq\bigcap_{F\in\Z_{\geq1}(\M)}\Delta_F$, which is to say that the
  lattice of cyclic flats has a unique minimal element of positive rank.
  For Schubert matroids, this is the case because $\Z(\M_S)$ is a chain.
\end{proof}

\begin{question}
  Are there other families of
  matroids for which $\Z_{\geq1}(\M)$ has a unique minimal element,
  and that satisfy \eqref{eq:res formula}?  
\end{question}

\section{Further directions}\label{sec:last}
We showed that Schubert matroids satisfy the convolution resonance formula
\eqref{eq:res formula}, and it seems natural to ask which other matroids
do as well.

For rank-$3$ matroids, the question has a straightforward answer.
We refer to the paper of Falk and Yuzvinsky~\cite{FY07} for the notion of
multinets.  For further discussion, see \cite[\S5]{De16}.  We will not
recall the definition here but simply record the relevant property.  They
show that the components of $\R^1(\M)$ are all linear subspaces.  These include
$\Delta_F$ for each connected flat $F$ (the \emph{local components}),
as well as some additional
components indexed by multinets, the existence of which imposes strong
constraints on $\M$. 
\begin{theorem}[Thm.\ 2.6, 2.7, \cite{FY07}]
  If $\M$ is a simple, rank-$3$ $\CC$-realizable matroid,
  \[
  \R^1(\M)\supseteq \bigcup_{F\in\L^\irr_2(\M)}\Delta_F
  \]
  is a strict containment if and only if $\M$ admits a multinet.
\end{theorem}
\begin{corollary}
  A simple, rank-$3$ matroid $\M$ realizable over $\CC$ satisfies the convolution
  resonance formula if and only if $\M$ admits no multinet.
\end{corollary}
\begin{proof}
  For nonzero $w\in \Delta_F$, we note $H^0(\os(\M),w)=0$ because $\M$ is
  simple.  Then we use \cite[Prop.\ 4.1]{De16} to reduce to
  a rank-$2$ computation, 
  $H^1(\os(\M),w)=H^1(\os(\M|_F),w)\cong \k^{\beta(\M|_F)}$.
  Noting that the Euler characteristic of $(\os(\M),w)$ is (independently of
  $w$) equal to $\beta(\M)$ determines $H^2(\os(\M),w)$, which can be seen
  to agree with the contributions from $F$ and $E$ in \eqref{eq:res formula}.

  Conversely, if $w$ belongs to a component of $\R^1(\M)$ indexed by a multinet,
  then $w\not\in \Delta_F$ for any rank-$2$ connected flat $F$.
  In that case, the right-hand side of \eqref{eq:res formula} is zero in
  degree $1$, hence not equal to $H^1(\os(\M),w)$.
  \end{proof}
\begin{remark}
  The $3$-whirl is easily seen not to admit a multinet, so it satisfies
  the convolution resonance formula; however, it is not arboreal
  (Example~\ref{ex:arboreal}).
\end{remark}

Not every arboreal matroid satisfies the convolution resonance formula
\eqref{eq:res formula}, however.  Since we have shown the Orlik--Solomon
algebra is monomial (Theorem~\ref{thm:degeneration}), the task remains to
interpret the Hochster-type resonance formula of \cite{PS09} in terms of 
matroid data such as cyclic flats.  We leave this for future work.
\begin{example}
  Let $\M$ be the rank-$4$ truncation of $\U_{23}\oplus\U_{23}\oplus\U_{11}$
  on the ground set $[6]$ and $w=w_1+w_2$, where $w_1\in\Delta_{[0,2]}$
  and $w_2\in\Delta_{[3,5]}$ are nonzero.
  A Macaulay2~\cite{M2} check shows $H^2(\os(\M),w)\neq0$, while
  the right-hand side of \eqref{eq:res formula} is zero in degree $2$.
\end{example}

\subsection*{Acknowledgements}
The first author thanks NSERC of Canada for support.
The authors acknowledge the use of Claude Fable 5.1 for help with the figures
and proofreading.

\bibliographystyle{amsalpha}
\bibliography{schres}
\end{document}

%% file: sch_arxiv.fig1.tex
\[
\begin{tikzpicture}[baseline=(current bounding box.south), scale=0.625, line cap=round, line join=round]
  \fill[blue!10] (0,0) -- (0,1) -- (1,1) -- (1,2) -- (1,3) -- (2,3) -- (2,4) -- (3,4) -- (3,0) -- cycle;
  \draw[blue, line width=1pt] (0,0) -- (0,1) -- (1,1) -- (1,2) -- (1,3) -- (2,3) -- (2,4) -- (3,4);
  \fill[red] (0,0) circle (2pt);
  \fill[red] (3,4) circle (2pt);
  \fill[orange!10] (3,0) -- (3,3) -- (4,3) -- (4,2) -- (5,2) -- (5,1) -- (6,1) -- (6,0) -- cycle;
  \draw[step=1, gray!60, thin] (0,0) grid (7,4);
  \draw[gray!85, semithick] (0,0) -- (0,1) -- (1,1) -- (1,2) -- (1,3) -- (2,3) -- (2,4) -- (3,4);
  \draw[orange!75!black, semithick] (3,4) -- (3,3) -- (4,3) -- (4,2) -- (5,2) -- (5,1) -- (6,1) -- (6,0) -- (7,0);
  \draw[red, line width=1pt] (0,0) -- (1,0) -- (1,1) -- (2,1) -- (2,2) -- (2,3) -- (3,3) -- (4,3);
  \fill[red] (0,0) circle (2pt);
  \fill[red] (4,3) circle (2pt);
  \node[below left, font=\tiny] at (0,0) {$(0,0)$};
  \node[font=\tiny] (endlab) at (4.9,4.4) {$(4,3)$};
  \draw[-{Stealth[length=3.5pt,width=3pt]}, gray!45!black, thin, shorten >=6pt] (5.05,3.9) to[bend left=25] (4,3);
  \node[font=\tiny] (endlab) at (3.5,4.4) {$(3,4)$};
  \node[left, red, font=\tiny] at (1,0.5) {2};
  \node[left, red, font=\tiny] at (2,1.5) {4};
  \node[left, red, font=\tiny] at (2,2.5) {5};
  \node[left, gray!40!black, font=\tiny] at (0,0.5) {1};
  \node[left, gray!40!black, font=\tiny] at (1,1.5) {3};
  \node[left, gray!40!black, font=\tiny] at (1,2.5) {4};
  \node[left, gray!40!black, font=\tiny] at (2,3.5) {6};
\end{tikzpicture}%
\hspace{10mm}%
\begin{tikzpicture}[baseline=(current bounding box.south), scale=0.625, line cap=round, line join=round]
  \fill[blue!12] (0,0) -- (0,1) -- (1,1) -- (1,2) -- (1,3) -- (2,3) -- (2,4) -- (3,4) -- (3,0) -- cycle;
  \draw[blue, line width=1pt] (0,0) -- (0,1) -- (1,1) -- (1,2) -- (1,3) -- (2,3) -- (2,4) -- (3,4);
  \fill[orange!20] (3,0) -- (3,3) -- (4,3) -- (4,2) -- (5,2) -- (5,1) -- (6,1) -- (6,0) -- cycle;
  \draw[step=1, gray!60, thin] (0,0) grid (7,4);
  \draw[gray!85, semithick] (0,0) -- (0,1) -- (1,1) -- (1,2) -- (1,3) -- (2,3) -- (2,4) -- (3,4);
  \draw[orange!75!black, semithick] (3,4) -- (3,3) -- (4,3) -- (4,2) -- (5,2) -- (5,1) -- (6,1) -- (6,0) -- (7,0);
  \draw[red, dashed, line width=1pt] (1,1) -- (2,1) -- (2,2) -- (2,3) -- (3,3) -- (4,3);
  \draw[red, fill=white, line width=0.8pt] (4,3) circle (1.7pt);
  \node[left, red!50, font=\tiny] at (2,1.5) {4};
  \node[left, red!50, font=\tiny] at (2,2.5) {5};
  \draw[-Implies, double, line width = .75pt]
       (1.50,1.18) -- (1.50,1.25) arc[start angle=0, end angle=90, radius=0.21] -- (1.18,1.50);
  \draw[red, line width=1pt] (0,0) -- (1,0) -- (1,1) -- (1,2) -- (1,3) -- (1,4) -- (2,4) -- (3,4);
  \fill[red] (0,0) circle (2pt);
  \fill[red] (3,4) circle (2pt);
\node[cross out, draw, minimum size=1mm, inner sep=0pt] (Xnode) at (1,3) {}; 
  \node[below left, font=\tiny] at (0,0) {$(0,0)$};
  \node[left, red, font=\tiny] at (1,0.5) {2};
  \node[left, red, font=\tiny] at (1,1.5) {3};
  \node[left, red, font=\tiny] at (1,2.5) {4};
  \node[left, red, font=\tiny] at (1,3.5) {5};
  \node[left, gray!40!black, font=\tiny] at (0,0.5) {1};
  \node[left, gray!40!black, font=\tiny] at (2,3.5) {6};
    \node[font=\tiny] (endlab) at (3.5,4.4) {$(3,4)$};
\end{tikzpicture}
\]

%% file: sch_arxiv.fig3.tex
\tikzset{
  bdry/.style    ={black, line width=1.4pt, line cap=round, line join=round},
  cont/.style    ={black, line width=2pt,   line cap=round, dash pattern=on 0pt off 7pt},
  bdots/.style   ={black, line width=2pt,   line cap=round, dash pattern=on 0pt off 7pt},
  ipath/.style   ={red,   line width=1.4pt, line cap=round},
  iedge/.style   ={red,   line width=1.4pt, line cap=round, line join=round},
  reddots/.style ={red,   line width=2pt,   line cap=round, dash pattern=on 0pt off 7pt},
  bluedots/.style={blue,  line width=2pt,   line cap=round, dash pattern=on 0pt off 7pt},
  bluep/.style   ={blue,  line width=1.3pt, line cap=round, line join=round, dash pattern=on 5pt off 4pt},
  i1/.style      ={blue,  line width=1.4pt, line cap=round},
}
\centering
\begin{subfigure}[b]{0.32\textwidth}
  \centering
  \begin{tikzpicture}[baseline={(0,0)}, scale=0.6]
    \fill[blue!12] (-1.5,-1.5) -- (1.5,-1.5) -- (1.5,1.5) -- (0,1.5) -- (0,0) -- (-1.5,0) -- cycle;
            \draw[step=1.5, gray!60, thin] (-1.5,-1.5) grid (1.5,1.5);
    \draw[cont, line width=1.2pt] (-1.65,-0.13) -- (-2.65,-1.13);
    \draw[cont, line width=1.2pt] (0.13,1.63) -- (1.13,2.63);
    \draw[ipath, line width=.75pt] (-1.5,0) -- (0,0);
    \draw[i1, line width=.75pt] (0,0) -- (0,1.5);
    \draw[black, fill=white, line width=.8pt] (0,0) circle (3pt);
    \node[black, font=\scriptsize] at (-0.75,-0.4) {$i$};
    \node[black, font=\scriptsize] at (0.65,0.85)   {$i{+}1$};
      \end{tikzpicture}
  \caption{$i$ horizontal, $i+1$ vertical.\label{fig:case1}}
  \label{fig:corner-both}
\end{subfigure}
\hfill
\begin{subfigure}[b]{0.32\textwidth}
  \centering
  \begin{tikzpicture}[baseline={(0,0)}, scale=0.6]
    \fill[blue!12] (-1.5,-1.5) -- (1.5,-1.5) -- (1.5,1.5) -- (0,1.5) -- (0,0) -- (-1.5,0) -- cycle;
    \draw[step=1.5, gray!60, thin] (-1.5,-1.5) grid (1.5,1.5);
    \draw[cont, line width=1.2pt] (-1.65,-0.13) -- (-2.65,-1.13);
    \draw[cont, line width=1.2pt] (0.13,1.63) -- (1.13,2.63);
    \draw[bdry, line width=.75pt] (-1.5,0) -- (0,0);
    \draw[i1, line width=.75pt] (0,0) -- (0,1.5);
    \draw[ipath, line width=.75pt] (0,0) -- (0,-1.5);
    \draw[black, fill=white, line width=.8pt] (0,0) circle (3pt);
    \node[black, font=\scriptsize] at (0.65,0.85)  {$i{+}1$};
    \node[red,   font=\scriptsize] at (0.38,-0.8) {$i$};
  \end{tikzpicture}
  \caption{$i$ and $i+1$ vertical.\label{fig:case2}}
  \label{fig:corner-cases}
\end{subfigure}
\hfill
\begin{subfigure}[b]{0.32\textwidth}
  \centering
  \begin{tikzpicture}[baseline={(2,1.25)}, scale=0.6]
    \fill[blue!12] (1,0) -- (3,0) -- (3,5) -- (2,5) -- (2,4) -- (1,4) -- cycle;
            \draw[step=1, gray!60, thin] (1,0) grid (3,5);
    \draw[bdry, line width=.75pt] (1,4) -- (2,4) -- (2,5);
    \draw[bdots, line width=1.2pt] (0.82,3.82) -- (0.1,3.1);
    \draw[bdots, line width=1.2pt] (2.18,5.18) -- (2.9,5.9);
    \draw[iedge, line width=.75pt] (1,0) -- (2,0) -- (2,4);
    \draw[reddots, line width=1.2pt] (2.4,1.7) -- (2.4,2.7);
    \draw[i1, line width=.75pt] (2,4) -- (3,4);
    \draw[bluep, line width=.75pt] (1,0) -- (1,5) -- (2,5);
    \draw[bluedots, line width=1.2pt] (0.65,2.13) -- (0.65,3.13);
  \draw[-Implies, double, line width = .75pt]
         (1.62,0.14) -- (1.62,0.32) arc[start angle=0, end angle=90, radius=0.2] -- (1.18,0.5);
  \draw[-Implies, double, line width = .75pt]
         (2.62,4.14) -- (2.62,4.32) arc[start angle=0, end angle=90, radius=0.2] -- (2.18,4.5);
    \draw[black, fill=white, line width=.75pt] (2,4) circle (3pt);
    \node[red,  font=\scriptsize] at (2.42,3.5)  {$i$};
    \node[blue, font=\scriptsize] at (3.5,4.4)   {$i{+}1$};
    \node[red,  font=\scriptsize] at (2.70,0.55) {$j{+}1$};
    \node[red,  font=\scriptsize] at (1.5,-0.5)  {$j$};
    \node[blue, font=\scriptsize] at (0.55,0.55) {$j$};
    \node[blue, font=\scriptsize] at (0.35,1.55) {$j{+}1$};
    \node[blue, font=\scriptsize] at (0.6,4.4)   {$i$};
  \end{tikzpicture}
  \caption{$i$ vertical, $i+1$ horizontal.\label{fig:case3}}
  \label{fig:exchange}
\end{subfigure}

%% file: sch_arxiv.fig2.tex
\[
\begin{tikzpicture}[scale=0.625, line cap=round, line join=round]
  \fill[blue!12] (0,0) -- (0,1) -- (1,1) -- (1,2) -- (1,3) -- (2,3) -- (2,4) -- (3,4) -- (3,0) -- cycle;
  \fill[orange!20] (3,0) -- (3,3) -- (4,3) -- (4,2) -- (5,2) -- (5,1) -- (6,1) -- (6,0) -- cycle;
  \draw[step=1, gray!60, thin] (0,0) grid (7,4);
  \draw[gray!85, semithick] (0,0) -- (0,1) -- (1,1) -- (1,2) -- (1,3) -- (2,3) -- (2,4) -- (3,4);
  \draw[orange!75!black, semithick] (3,4) -- (3,3) -- (4,3) -- (4,2) -- (5,2) -- (5,1) -- (6,1) -- (6,0) -- (7,0);
  \draw[red, line width=1pt] (0,0) -- (1,0) -- (1,1) -- (2,1) -- (2,2) -- (2,3);
  \draw[red, dashed, line width=.75pt] (0,0) -- (0,1) -- (1,1) -- (1,2) -- (1,3) -- (2,3);
  \draw[blue, line width=1pt] (2,3) -- (2,4) -- (3,4);
  \fill[red] (0,0) circle (2pt);
  \fill[blue] (3,4) circle (2pt);
  \draw[red, fill=white, line width=0.8pt] (2,3) circle (1.7pt);
  \draw[-Implies, double, line width = .75pt, shorten >=5pt, shorten <=5pt]
 (1,2) -- (2,2);
  \node[left, red!50, font=\tiny] at (0,0.5) {1};
  \node[left, red!50, font=\tiny] at (1,1.5) {3};
  \node[left, red!50, font=\tiny] at (1,2.5) {4};
  \node[left, red, font=\tiny] at (1,0.5) {2};
  \node[left, red, font=\tiny] at (2,1.5) {4};
  \node[left, red, font=\tiny] at (2,2.5) {5};
  \node[left, blue, font=\tiny] at (2,3.5) {6};
\end{tikzpicture}
\]